\documentclass{amsart}
\usepackage{amsmath,amsfonts,amssymb,amscd,amsthm,amsbsy,epsf}
\newtheorem{thm}{Theorem} 
\newtheorem{cor}{Corollary} 
\newtheorem*{Cor}{Corollary} 

\newtheorem*{lem}{Lemma}
\theoremstyle{remark}
\newtheorem*{rem}{Remark}
\def\F{{\mathcal F}}
\def\S{{\mathcal S}}
\def\real{{\mathbb  R}}

\def\hatf{\widehat f}
\def\ds{\displaystyle}
\def\step#1{\medskip\noindent$\underline{\text{Step #1}}$.\enspace}

\def\smallnegint{\mathop{\int\mkern-13mu
        \raise.5ex\hbox{${\scriptscriptstyle\diagup}$}}\nolimits}
\def\too{\longrightarrow}
\begin{document}

\title{Weighted inequalities and Stein-Weiss potentials}
\author{William Beckner}
\address{Department of Mathematics, The University of Texas at Austin,
1 University Station C1200, Austin TX 78712-0257 USA}
\email{beckner@math.utexas.edu}
\begin{abstract}
Sharp extensions of Pitt's inequality and bounds for Stein-Weiss fractional 
integrals are obtained that incorporate gradient forms and vector-valued
operators.
Such results include Hardy-Rellich inequalities.
\end{abstract}
\maketitle

Weighted inequalities provide quantitative information to characterize 
integrability for differential and integral operators and intrinsically are 
determined by their dilation character.
In the classical context, weighted inequalities for the Fourier transform 
provide a natural measure of uncertainty.
For functions on $\real^n$ the issue is the balance between the relative 
size of a function and its Fourier transform at infinity.
An inequality that illustrates this principle at the spectral level is 
Pitt's inequality:
\begin{equation}\label{eq:spectral}
\int_{\real^n} \Phi (1/|x|) |f(x)|^2\,dx 
\le C_\Phi \int_{\real^n} \Phi (|y|) |\hatf (y)|^2\,dy
\end{equation}
where $\Phi$ is an increasing function, the function $f$ is in the 
Schwartz class $\S(\real^n)$ and the Fourier transform is defined by 
\begin{equation*}
(\F f)(y) = \hatf (y) = \int_{\real^n} e^{2\pi ixy} f(x)\,dx\ .
\end{equation*}
Such inequalities may be fully determined by dilation invariance, and some
cases may be realized with explicit gradient forms as Hardy-Rellich
inequalities. 
In earlier work (see \cite{3}) the effective calculation for the constant 
in Pitt's inequality was reduced to Young's inequality for convolution on 
a non-compact unimodular group. 
The objective here will be to study more general forms of Pitt's inequality
\begin{equation}\label{general-Pitt}
\int_{\real^n} \Phi (1/|x|) |\nabla f|^2\,dx
\le 4\pi^2 D_\Phi \int_{\real^n} \Phi (|y|) |y|^2 |\hatf (y)|^2\,dy
\end{equation}
using the structure of Stein-Weiss potentials and convolution estimates 
in concert with the Hecke-Bochner representation for $L^2 (\real^n)$. 
The previous work is described by the following three theorems.

\begin{thm}[Pitt's inequality]\label{thm:Pitt}
For $f\in \S(\real^n)$ and $0\le \alpha <n$
\begin{gather}\label{eq:pitt}
\int_{\real^n} |x|^{-\alpha} |f(x)|^2\,dx 
\le C_\alpha \int_{\real^n} |y|^\alpha |\hatf (y)|^2\,dy\\
\noalign{\vskip6pt}
C_\alpha = \pi^\alpha \Big[  \Gamma \Big(\frac{n-\alpha}4\Big) \Big/
\Gamma \Big(\frac{n+\alpha}4\Big)\Big]^2\ .\notag
\end{gather}
\end{thm}

\noindent
Since the above inequality becomes an identity for $\alpha=0$, a 
differentiation argument provides a logarithmic form that controls the 
uncertainty principle by using dimensional asymptotics.

\begin{thm}[logarithmic uncertainty]\label{thm:log}
For $f\in \S(\real^n)$ 
\begin{gather}\label{eq:log}
\int_{\real^n} \ln |x|\, |f(x)|^2\,dx 
+ \int_{\real^n} \ln |y|\, |\hatf (y)|^2\,dy 
\ge D\int |f(x)|^2\,dx \\
\noalign{\vskip6pt}
D= \psi (n/4) - \ln \pi\ ,\qquad 
\psi (t) = \frac{d}{dt} \ln \Gamma (t)\ .\notag
\end{gather}
\end{thm}

\noindent 
Logarithmic integrals are indeterminate so $D$ may take negative values.
The proof of Pitt's inequality \eqref{eq:pitt} follows from a sharp 
estimate for an equivalent integral realization as a Stein-Weiss 
fractional integral on $\real^n$.

\begin{thm}\label{thm:pittsproof}
For $f\in L^2 (\real^n)$ and $0<\alpha <n$
\begin{gather}\label{eq:pittsproof}
\Big|\int_{\real^n\times\real^n} \mkern-18mu 
f(x) \frac1{|x|^{\alpha/2}} 
\frac1{|x-y|^{n-\alpha}} \frac1{|y|^{\alpha/2}} f(y)\,dx\,dy\Big|
\le B_\alpha \int_{\real^n} |f(x)|^2\,dx\\
\noalign{\vskip6pt}
B_\alpha = \pi^{n/2} 
\Big[\Gamma\Big(\frac{\alpha}2\Big)\Big/ 
\Gamma\Big(\frac{n-\alpha}2\Big)\Big]
\Big[\Gamma\Big(\frac{n-\alpha}4\Big)\Big/
\Gamma\Big(\frac{n+\alpha}4\Big)\Big]^2\ .\notag
\end{gather}
\end{thm}

\noindent
The development of the sharp estimate for the Stein-Weiss integral rests 
on using symmetrization to reduce the problem to radial functions, and 
then as a consequence of dilation invariance the estimate can be converted 
to Young's inequality for convolution on the multiplicative group 
$\real_+$ (or alternatively on $\real$).
$$\| \varphi * f\|_{L^2(G)} \le \|\varphi\|_{L^1(G)} \|f\|_{L^2(G)}$$
When $\varphi$ is non-negative, the inequality is sharp with no extremal 
functions.
One objective here will be to extend this inequality to cases when 
$\psi$ takes both positive and negative values by using the 
Hecke-Bochner formulas.

Since the form 
$$\int_{\real^n} |\hatf(y)|^2 |y|^\alpha\,dy 
= (2\pi)^{-\alpha} \int_{\real^n} |(-\Delta)^{\alpha/4}f|^2\,dx$$
can be regarded in the family of gradient estimates, Pitt's inequality 
has been characterized as a Hardy-Rellich inequality in some parts 
of the recent literature with alternative proofs and extensions 
(see \cite1, \cite7, \cite{14}, and \cite{16}). 
Some natural generalizations of Pitt's inequality can now be viewed in 
the context of the arguments developed in \cite3 where proofs of 
Theorems~\ref{thm:Pitt}, \ref{thm:log} and \ref{thm:pittsproof} are given.

\section{Pitt's inequality with gradient terms}

\begin{thm}\label{thm:gradient}         
For $f\in \S(\real^n)$ and $0<\alpha <n$, $n>1$
\begin{gather}\label{eq:gradient}
\int_{\real^n} |\nabla f|^2|x|^{-\alpha}\,dx
\le 4\pi^2 D_\alpha \int_{\real^n} |\hatf(y)|^2 |y|^{\alpha+2}\,dy\\
\noalign{\vskip6pt}
D_\alpha = \pi^\alpha \max_k 
\bigg\{\Big[ \Gamma\Big(\frac{n+2k-\alpha+2}4\Big) \Big/
\Gamma\Big( \frac{n+2k+\alpha +2}4\Big)\Big]^2 
\Big( 1+\frac{4k\alpha}{(n+2k-\alpha -2)^2}\Big)\bigg\}\notag
\end{gather}
\end{thm}

\noindent 
By convention the last term in the line above is one when $k=0$. 
This result is interesting for several aspects:
\begin{itemize}
\item[(a)] for the gradient term with radial functions (i.e. $k=0$) 
the constant is reduced from that in equality~\eqref{eq:pitt} since 
$\Gamma (x+\beta)/\Gamma (y+\beta)$ is decreasing in $\beta$ for $x<y$;
this reduction in constant is also apparent by the Plancherel theorem;
\item[(b)] the corresponding Stein-Weiss integral does not have a positive 
kernel though it does have symmetry in the angular variables;
\item[(c)] the limiting logarithmic uncertainty is sharper 
\begin{gather*}
\int_{\real^n} \ln |x| \, |\nabla f|^2\,dx 
+ \int_{\real^n} \ln |y| (4\pi^2 |y|^2) |\hatf (y)|^2\,dy 
\ge E\int_{\real^n} |\nabla f|^2\,dx\\
\noalign{\vskip6pt}
E = \begin{cases}
\ds \psi\Big(\frac32\Big) - \ln \pi -1\ ,&n=2\\
\noalign{\vskip6pt}
\ds \psi\Big(\frac{n}{4} +\frac12\Big) - \ln\pi\ ,&n\ge3
\end{cases}
\end{gather*}
\item[(d)] for $\alpha=2$ and $n>4$ one obtains the Hardy-Rellich inequality 
\begin{equation}\label{eq:hardy-rellich}
\int_{\real^n} |\nabla f|^2 |x|^{-2} \,dx 
\le \frac4{n^2} \int_{\real^n} |\Delta f|^2\,dx\ ;
\end{equation}
\item[and (e)] 
for some values of the parameters $n$ and $\alpha$, the sharp bound is 
obtained by considering non-radial functions.
A sharper version of \eqref{eq:hardy-rellich} appears in \cite{14}.
\end{itemize}

By using the Fourier transform on Riesz potentials 
$$\F\big[|x|^{-\alpha}\big] = \pi^{-\frac{n}{2}\, +\alpha}
\Big[\Gamma\Big(\frac{n-\alpha}2\Big)\Big/ \Gamma\Big(\frac{\alpha}2\Big)\Big]
|x|^{-n+\alpha}\ , $$
one easily sees that Pitt's inequality in Theorem~\ref{thm:gradient} 
is equivalent to a Stein-Weiss fractional integral inequality on $\real^n$.

\begin{thm}\label{thm5}
For $f\in L^2 (\real^n)$ and $0<\alpha <n$, $n>1$
\begin{gather}\label{eq8}
\Big| \int_{\real^n\times\real^n}\mkern-18mu  f(x)
\frac1{|x|^{\alpha/2}} \ 
\frac{x\cdot y}{|x|\,|y|}\  
\frac 1{|x-y|^{n-\alpha}} \ 
\frac1{|y|^{\alpha/2}} 
f(y) \,dx\,dy \Big|\\
\noalign{\vskip6pt}
\le \Big[ \pi^{\frac{n}2 -\alpha} \Gamma\Big(\frac{\alpha}2\Big)\Big/
\Gamma\Big(\frac{n-\alpha}2\Big)\Big]\, D_\alpha        \notag
\int_{\real^n} |f|^2\,dx 
\end{gather}
with $D_\alpha$ as in Theorem~\ref{thm:gradient}.
\end{thm}

\begin{proof}[Proof of Theorems~\ref{thm:gradient} and \ref{thm5}] 
The main difficulty with showing these estimates is that the Stein-Weiss 
kernel is not positive. 
This obstruction is addressed in two steps.

\step{1} 
Assume $f$ is radial and set $t= |x|$, $h(t) = |x|^{n/2} f(x)$; 
then \eqref{eq8} reduces to 
\begin{equation}\label{eq9}
\Big|\int_{\real_+\times\real_+} \mkern-18mu
h(t)\psi_\alpha(s/t)h(s)\frac{ds}s\, 
\frac{ds}t\Big| \le D_\alpha 
\Big[\Gamma\Big(\frac{n}2\Big)
\Gamma\Big(\frac{\alpha}2\Big)\Big/ 2\pi^\alpha 
\Gamma\Big(\frac{n-\alpha}2\Big)\Big] 
\int_{\real_+} |h|^2 \frac{dt}t
\end{equation}
with 
$$\psi_\alpha(t) = \int_{S^{n-1}} \xi_1 
\Big[ t+\frac1t - 2\xi_1\Big]^{-(n-\alpha)/2} \,d\xi$$
where $d\xi$ denotes normalized surface measure, $\xi_1$ is the first 
component of $\xi$ and $D_\alpha$ is specified to be the best constant 
in \eqref{eq:gradient}. 
Note that by monotonicity and symmetry $\psi_\alpha (t)$ is positive, 
and moreover any kernel derived from \eqref{eq:gradient} will be 
positive-definite.
Then using Young's inequality
$$\|\psi_\alpha * h\|_{L^2(\real_+)} 
\le \|\psi_\alpha\|_{L^1(\real_+)} \|h\|_{L^2 (\real_+)}$$ 
which is sharp, it suffices to calculate the $L^1$ norm of $\psi_\alpha$ 
to obtain the constant $D_\alpha$
$$\|\psi_\alpha\|_{L^1(\real_+)} = D_\alpha 
\Big[\Gamma \Big(\frac{n}2\Big)\Gamma(\alpha/2)
\Big/2\pi^\alpha \Gamma\Big(\frac{n-\alpha}2\Big)\Big]\ .$$
To compute this integral, observe that 
\begin{equation*}
\begin{split}
\|\psi_\alpha\|_{L^1(\real_+)} 
& = \int_0^\infty \Big[ \int_{S^{n-1}} \xi_1 
\Big[ t+\frac1t - 2\xi_1\Big]^{-(n-\alpha)/2}\,d\xi\Big]\, \frac{dt}t\\
\noalign{\vskip6pt}
& = \Big[ \frac{2\pi^{n/2}}{\Gamma(n/2)}\Big]^{-1} 
\int_{\real^n} \frac{x\cdot y}{|x|\, |y|}\ 
\frac1{|x-y|^{n-\alpha}}\ 
\frac1{|y|^{(n+\alpha)/2}}\, dy
\end{split}
\end{equation*}
for $|x|=1$ and $n>1$. 
The second integral will be calculated for the set of values $n-2>\alpha>0$. 
But notice that the first integral is an analytic function of the parameters
$n\ge2$ and $\beta = n-\alpha$ for some range of values.
Hence any computation for some parameter domain will determine by 
analytic continuation the value of $\|\psi_\alpha\|_{L^1(\real_+)}$ 
for the desired parameter interval $0<\alpha <n$ with $n>1$. 
Then (noting that $|x|=1$) 
\begin{equation*}
\begin{split} 
&\int_{\real^n} \frac{2x\cdot y}{|x|\, |y|}\ |x-y|^{-(n-\alpha)} 
|y|^{-(n+\alpha)/2}\,dy  = \\
\noalign{\vskip6pt}
&\qquad  \int_{\real^n} |x-y|^{-(n-\alpha)} |y|^{-(n+\alpha)/2\ -1}\,dy 
+ \int |x-y|^{-(n-\alpha)} |y|^{-(n+\alpha)/2\ +1}\,dy\\
\noalign{\vskip6pt}
&\qquad - \int_{\real^n} |x-y|^{-(n-\alpha -2)} 
|y|^{-(n+\alpha)/2-1}\,d y 
= I_1 + I_2 - I_3
\end{split}
\end{equation*}
These integrals are computed by using the formula for the convolution 
of two Riesz potentials
\begin{equation} \label{eq10}
|x|^{-\beta} * |x|^{-\delta} = \pi^{n/2} 
\left[ \frac{\Gamma(\frac{n-\beta}2) \Gamma(\frac{n-\delta}2) 
\Gamma(\frac{\beta+\delta -n}2)}
{\Gamma (\frac{\beta}2) \Gamma (\frac{\delta}2) 
\Gamma(\frac{2n-\beta-\delta}2)}
\right] |x|^{-(\beta+\delta -n)}
\end{equation}
with $0<\beta <n$, $0<\delta<n$ and $n<\beta +\delta <2n$. 
Then 
\begin{equation*}
\begin{split}
I_1 & = \pi^{n/2} \left[ 
\frac{\Gamma (\frac{\alpha}2) \Gamma(\frac{n-\alpha-2}4) 
\Gamma (\frac{n-\alpha+2}4) }
{\Gamma (\frac{n-\alpha}2) \Gamma (\frac{n+\alpha+2}4) 
\Gamma (\frac{n+\alpha-2}4)}\right]
\\
\noalign{\vskip6pt}
I_2 & = \pi^{n/2} \left[ 
\frac{\Gamma (\frac{\alpha}2) \Gamma(\frac{n-\alpha+2}4) 
\Gamma (\frac{n-\alpha-2}4)}
{\Gamma (\frac{n-\alpha}2) \Gamma (\frac{n+\alpha-2}4) 
\Gamma (\frac{n+\alpha+2}4)}\right]
\\
\noalign{\vskip6pt}
I_3 & = \pi^{n/2} \left[ 
\frac{\Gamma (\frac{\alpha+2}2) \Gamma(\frac{n-\alpha-2}4) 
\Gamma (\frac{n-\alpha-2}4)}
{\Gamma (\frac{n-\alpha-2}2) \Gamma (\frac{n+\alpha+2}4) 
\Gamma (\frac{n+\alpha+2}4) }\right]
\end{split}
\end{equation*}
and 
\begin{equation*}
\begin{split}
I_1 + I_2 - I_3 
& = \pi^{n/2} 
\frac{\Gamma (\frac{\alpha}2)}{\Gamma (\frac{n-\alpha}2)}
\left[ \frac{\Gamma(\frac{n-\alpha -2}4)}
{\Gamma (\frac{n+\alpha +2}4)}\right]^2 
\left[ 2\Big(\frac{n-\alpha-2}4\Big) 
\Big(\frac{n+\alpha-2}4\Big) 
- \frac{\alpha}2 \Big(\frac{n-\alpha -2}2\Big) \right]\\
\noalign{\vskip6pt}
& = \pi^{n/2} 
\frac{\Gamma (\frac{\alpha}2)}{\Gamma (\frac{n-\alpha}2)}
\left[ \frac{\Gamma(\frac{n-\alpha -2}4)}
{\Gamma (\frac{n+\alpha +2}4)}\right]^2 
\left[ 2\Big(\frac{n-\alpha-2}4\Big)^2\right] 
= 2\pi^{n/2} 
\frac{\Gamma (\frac{\alpha}2)}{\Gamma (\frac{n-\alpha}2)}
\left[ \frac{\Gamma (\frac{n-\alpha+2}4)}{\Gamma(\frac{n+\alpha+2}4)}
\right]^2 .
\end{split}
\end{equation*}
This demonstrates that 
$$\|\psi_\alpha\|_{L^1 (\real_+)} = 
\frac{\Gamma (\frac{n}2) \Gamma(\frac{\alpha}2)}
{2\Gamma (\frac{n-\alpha}2)} 
\left[ \frac{\Gamma (\frac{n-\alpha+2}4)}{\Gamma(\frac{n+\alpha+2}4)}
\right]^2$$
and that for radial functions  in \eqref{eq:gradient}
$$D_\alpha = \pi^\alpha \left[\frac{\Gamma (\frac{n-\alpha+2}4)}
{\Gamma(\frac{n+\alpha +2}4)}\right]^2 \ .$$

\step{2} 
The Hecke-Bochner representation for $L^2 (\real^n)$ is used to reduce the 
study of inequality \eqref{eq:gradient} 
$$\int_{\real^n} |\nabla f|^2 |x|^{-\alpha}\,dx 
\le 4\pi^2 D_\alpha \int_{\real^n} |\hatf (y)|^2 |y|^{\alpha +2}\,dy$$
to estimates for radial functions. 
For $f\in \S(\real^n)$
$$f(x) = \sum_{k=0}^\infty f_k (|x|) P_k(x)$$
where $P_k$ is a harmonic polynomial of degree $k$,
$$P_k (x) = |x|^k Y_k (\xi)\ ,\qquad 
\xi = \frac{x}{|x|}\ ,\qquad 
\int_{S^{n-1}} |Y_k (\xi)|^2\,d\xi = \frac{\omega_{n-1+2k}}{\omega_{n-1}}$$
$Y_k$ is a spherical harmonic of degree $k$, $\omega_m =$ surface area of 
the unit sphere $S^m$, and $d\xi$ is normalized surface measure on $S^{n-1}$. 
Then 
$$\int_{\real^n} |f|^2\,dx 
= \sum_{k=0}^\infty |f_k(|x|)|^2\,dx \ .$$
Let $\F_n$ denote the Fourier transform on $\real^n$. 
Bochner's relation for spherical harmonics is 
\begin{equation}\label{eq11} 
\F_n \big(f_k(|x|)P_k(x)\big) 
=  i^k \F_{n+2k}( f_k (|x|)) P_k\ .
\end{equation}
and the integral on the right-hand side of \eqref{eq:gradient} becomes 
$$ \int_{\real^n} |\hatf (y)|^2 |y|^{\alpha+2}\,dy 
= \sum_{k=0}^\infty \int_{\real^{n+2k}}|\hatf_k (|y|)|^2 
|y|^{\alpha +2}\,dy\ .$$
For the integral on the left-hand side in \eqref{eq:gradient}
\begin{equation*}
\begin{split}
\int_{\real^n} |\nabla f|^2 |x|^{-\alpha} \,dx 
& = \sum_{k = 0}^\infty \int_{\real^n} 
|\nabla [  f_k (|x|) P_k(x)]|^2 
|x|^{-\alpha} \,dx\\
\noalign{\vskip6pt}
&=\sum_{k=0}^\infty \int_{\real^{n+2k}} \mkern-24mu
|\nabla f_k(|x|)|^2 |x|^{-\alpha}\,dx
+ \sum_{k=1}^\infty k\alpha \int_{\real^{n+2k}}\mkern-24mu |f_k(|x|)|^2 
|x|^{-\alpha-2}\,dx\ .
\end{split}
\end{equation*}
Using inequality \eqref{eq:gradient} for radial functions from Step~1, 
$$\int_{\real^{n+2k}} \mkern-24mu 
|\nabla f_k (|x|)|^2 |x|^{-\alpha} \,dx 
\le 4\pi^{2+\alpha} 
\left[ \frac{\Gamma  (\frac{n+2k-\alpha +2}4)}
{\Gamma (\frac{n+2k+\alpha +2}4)}\right]^2 
\int_{\real^{n+2k}}\mkern-24mu 
|\hatf_k(y)|^2 |y|^{\alpha +2}\,dy$$
and using inequality \eqref{eq:pitt} 
$$
\int_{\real^{n+2k}} \mkern-24mu 
|f_k (|x|)|^2 |x|^{-\alpha-2} \,dx 
\le \pi^{2+\alpha} 
\left[ \frac{\Gamma  (\frac{n+2k-\alpha -2}4)}
{\Gamma (\frac{n+2k+\alpha +2}4)}\right]^2 
\int_{\real^{n+2k}}\mkern-24mu 
|\hatf_k(y)|^2 |y|^{\alpha +2}\,dy$$
one obtains 
\begin{equation*}
\begin{split}
&\int_{\real^n} 
|\nabla f|^2 |x|^{-\alpha} \,dx \\
&\qquad \le 4\pi^{2+\alpha}  \sum_{k=0}^\infty 
\left[\bigg[
\frac{\Gamma(\frac{n+2k-\alpha+2}4)}{\Gamma(\frac{n+2k+\alpha+2}4)}\bigg]^2
+ \frac{k\alpha}4 
\bigg[ \frac{\Gamma (\frac{n+2k-\alpha-2}4)}{\Gamma(\frac{n+2k+\alpha+2}4)}
\bigg]^2 \right] 
\int_{\real^{n+2k}}\mkern-24mu |\hatf_k(y)|^2 |y|^{\alpha+2}\,dy\\
&\qquad \le 4\pi^{2+\alpha} \max_k 
\left\{\bigg[
\frac{\Gamma (\frac{n+2k-\alpha+2}4)}{\Gamma(\frac{n+2k+\alpha+2}4)}
\bigg]^2 
\bigg[1+\frac{4k\alpha}{(n+2k-\alpha-2)^2}\bigg]\right\}
\int_{\real^n} |\hatf (y)|^2 |y|^{\alpha+2}\,dy
\end{split}
\end{equation*}
which demonstrates inequality \eqref{eq:gradient}, Theorem~\ref{thm:gradient}
and the equivalent Theorem~\ref{thm5}. 
\renewcommand{\qed}{}
\end{proof}

\begin{cor}\label{cor1}
For $\alpha=2$
\begin{equation}\label{eq12}
\frac{D_2}{\pi^2} = \begin{cases}
\ds \frac{144}{25}\ ,&\text{if $n=3$ ($k=1$ term)}\\
\noalign{\vskip6pt}
\ds \ \, \frac43\ ,&\text{if $n=4$ ($k=1$ term)}\\
\noalign{\vskip6pt}
\ds \ \frac{16}{n^2}\ ,&\text{if $n>4$ ($k=0$ term)}
\end{cases}
\end{equation}
\end{cor}

\noindent 
This result recovers a ``classical Hardy-Rellich inequality'' for $n>4$
(see remarks in \cite{14})
$$\int_{\real^n} |\nabla f|^2 |x|^{-2}\,dx 
\le \frac4{n^2} \int_{\real^n} |\Delta f|^2\,dx\ .$$
Notice that for $n=8$ this inequality is entirely elementary since 
integrating by parts 
$$\int |\nabla f|^2 |x|^{-2}\,dx 
= - \int f(\Delta f) |x|^{-2} \,dx
- 8\int |f|^2 |x|^{-4}\,dx
\le \frac1{16} \int |\Delta f|^2\,dx\ .$$

In comparing terms to evaluate constants explicitly, the following fact 
is useful:
for $0<x<y$, the ratio $\Gamma (x+\beta)/\Gamma (y+\beta)$ is decreasing 
for $\beta >0$. 
Set $F(\beta) = \ln \Gamma (x+\beta) - \ln \Gamma (y+\beta)$; then 
$$F' (\beta) = \psi (x+\beta)-\psi (y+\beta) 
= - (y-x) \sum_{k=0}^\infty (x+\beta +k)^{-1} 
(y+\beta +k)^{-1} <0\ .$$
For $n=2$, the expression 
$$\left[ \frac{\Gamma (\frac{n+2k-\alpha +2}4)}
{\Gamma (\frac{n+2k+\alpha +2}4)}\right]^2 
\left( 1+\frac{4k\alpha}{(n+2k-\alpha -2)^2}\right)$$
is decreasing for $k\ge1$ since both terms in the product are decreasing 
so the value of $D_\alpha$ is found by comparing the terms for $k=0$ and $k=1$.

\begin{cor}\label{cor2}
For $n=2$
\begin{equation}\label{eq13} 
D_\alpha = \pi^\alpha \left[ \frac{\Gamma (\frac32 -\frac{\alpha}4)}
{\Gamma (\frac32 +\frac{\alpha}4)}\right]^2 
\left( \frac{4+\alpha^2}{(2-\alpha)^2}\right)\ . 
\end{equation}
\end{cor}

\begin{proof} 
Set $\beta = \alpha/4$ with $0\le \beta <\frac12$ since $\alpha <n=2$. 
Then to see that the $k=1$ term is larger than the $k=0$ term, consider
the log of the ratio of these two terms and it suffices to show that
$$F(\beta) = \ln \Gamma \Big(\frac12-\beta\Big) +\ln\Gamma (1+\beta) 
- \ln\Gamma (1-\beta) - \ln \Gamma \Big(\frac32 +\beta\Big) 
+ \frac12 \ln \Big(\frac14+\beta^2\Big)$$
is positive for $\beta >0$ which will follow by showing that $F(\beta)$ 
is increasing with ${f(0)=0}$. 
\begin{equation*}
\begin{split}
F'(\beta) & = -4\Big(\frac12-\beta\Big) + \psi(1+\beta) + \psi(1-\beta) 
- \psi \Big(\frac32 +\beta\Big) + \frac{4\beta}{1+4\beta^2}\\
\noalign{\vskip6pt}
& =\Big(\frac12 +2\beta\Big) 
\sum_{k=0}^\infty \left[ \Big(\frac12 -\beta+k\Big)^{-1} 
(1+\beta+k)^{-1} - (1-\beta +k)^{-1} \Big(\frac32+\beta +k\Big)^{-1}\right] \\
\noalign{\vskip6pt}
&\qquad + \frac{4\beta}{1+4\beta^2} >0
\end{split}
\end{equation*}
$F'(\beta) >0$ and $F(0)=0$ ensure that $F(\beta)$ is positive.
\renewcommand{\qed}{}
\end{proof}

\begin{cor}\label{cor3}
For $n\ge 3$
\begin{equation}\label{eq14}
D_\alpha = \pi^\alpha \max_{k=0,1} 
\left\{ \bigg[\frac{\Gamma (\frac{n+2k-\alpha +2}4)}
{\Gamma (\frac{n+2k+\alpha +2}4)}\bigg]^2 
\left( 1+\frac{4k\alpha}{(n+2k-\alpha -2)^2}\right)\right\}\ .
\end{equation}
\end{cor}

\begin{proof}
The objective here is to show that only the $k=0$ and $k=1$ terms in 
Theorem~\ref{thm:gradient} need to be compared. 
That is, the only functions that are necessary to consider in 
Theorem~\ref{thm:gradient} are those contained in the span of functions 
with spherical harmonics up to degree one.
Set $\beta = 2k$ and show that 
$$G(\beta) = \ln \left\{\bigg[ \frac{\Gamma(\frac{n+\beta-\alpha+2}4)}
{\Gamma (\frac{n+\beta+\alpha+2}4)}\bigg] 
\left(1+\frac{2\beta\alpha}{(n+\beta -\alpha-2)^2}\right)^{1/2}\right\}$$
is decreasing for $\beta \ge2$ with $0<\alpha <n$; then
\begin{equation*}
\begin{split}
G'(\beta) & = \frac14 \left[ \psi\left(\frac{n+\beta-\alpha+2}4\right) 
-\psi \left(\frac{n+\beta+\alpha+2}4\right)\right]\\
\noalign{\vskip6pt}
&\qquad \quad+ \alpha \left[1-\frac{2\beta}{n+\beta-\alpha-2}\right]
\left[(n+\beta-\alpha -2)^2 + 2\beta\alpha\right]^{-1}\\
\noalign{\vskip6pt}
& = \alpha \Bigg[ -\frac18 \sum_{k=0}^\infty \bigg[ 
\left(\frac{n+\beta+2}4 +k\right)^2 - \frac{\alpha^2}{16}\bigg]^{-1} \\
\noalign{\vskip6pt}
&\qquad \quad
+ \frac{(n-\beta-\alpha-2)}{(n+\beta-\alpha-2)} 
\left[ (n+\beta-\alpha -2)^2 +2\beta\alpha\right]^{-1}\Bigg]\ .
\end{split}
\end{equation*}
This derivative is clearly negative for $n=3,4$ and $\beta \ge2$, 
$n>\alpha >0$ since $n-\beta-\alpha-2<0$.
Now consider $n\ge5$ with $0<\alpha <n$, and use a Riemann sum to 
approximate the first term from above:
\begin{gather*}
\frac18 \sum_{k=0}^\infty \left[\left(\frac{n+\beta+2}4 +k\right)^2
-\frac{\alpha^2}{16}\right]^{-1} 
> \frac18 \sum_{k=0}^\infty \left( \frac{n+\beta+2}4 +k\right)^{-2}\\
\noalign{\vskip6pt}
> \frac18 \int_0^\infty \left( \frac{n+\beta+2}4 +x\right)^{-2} dx 
= \frac12 (n+\beta +2)^{-1}\ .
\end{gather*}
Then 
$$G'(\beta) < \alpha \left[ -\frac12 (n+\beta+2)^{-1} 
+ \frac{n-\beta -\alpha-2}{n+\beta-\alpha-2} 
\left[ (n+\beta-2)^2 +\alpha^2 -2\alpha (n-2)\right]^{-1}\right]\ .$$ 
The right-hand expression is negative if 
$$\frac{n-\beta-\alpha-2}{n+\beta -\alpha-2} 
< \frac{(n+\beta-2)^2 +\alpha^2 -2\alpha (n-2)}{2(n+\beta+2)}\ ;$$
this is clearly the case if $n-\beta-\alpha-2<0$ so set 
$\delta = n-\beta -\alpha-2$ and consider the expression
$$\frac{-\delta}{\delta +2\beta} 
+ \frac{(2\beta +\delta+\alpha)^2 + \alpha^2 -2\alpha (\delta+\beta+\alpha)}
{2(\delta +2\beta +\alpha +4)}$$
or 
\begin{equation*}
\begin{split}
H(\delta) & = (\delta +2\beta) \left[(2\beta+\delta +\alpha)^2 
+\alpha^2 -2\alpha (\delta +\beta+\alpha)\right] 
- 2\delta (\delta +2\beta +\alpha +4)\\
\noalign{\vskip6pt}
& = (\delta+2\beta)(\delta +2\beta+\alpha)^2 -(\delta +2\beta) 
(\alpha^2 +2\alpha\delta + 2\alpha\beta) - 2\delta (\delta +2\beta +\alpha+4)\\
\noalign{\vskip6pt}
& = (\delta+2\beta) \left[\delta^2 +4\beta^2 + 2\beta\alpha +(4\beta-2\delta
\right] - 2\delta (\alpha +4) > 0
\end{split}
\end{equation*}
for $\beta\ge2$ and the positivity is clear for the case $\delta >0$. 
$H(\delta)>0$ implies that $G'(\beta) <0$, and this completes the proof 
of Corollary~\ref{cor3}.
\renewcommand{\qed}{}
\end{proof}

\begin{thm}[Hardy-Rellich trace inequality]\label{thm:trace}
For $f\in \S(\real^n)$, $n\ge2$
\begin{equation}\label{eq:trace1}
\int_{\real^n} |\nabla f|^2 |x|^{-1} \,dx 
\le \frac{D_1}{2\pi} \int_{\real^n} |(-\Delta)^{3/4} f|^2\,dx
\end{equation}
\begin{equation}\label{eq:trace2}
\frac{D_1}{2\pi} = \begin{cases} 
\ds \frac52 \left[\frac{\Gamma (\frac54)}{\Gamma (\frac74)}\right]^2\ ,&
\text{$n=2$ ($k=1$ term)}\\
\noalign{\vskip6pt}
\ds\frac{\pi}4\ ,&\text{$n=3$ ($k=1$ term)}\\
\noalign{\vskip6pt}
\ds\frac12 \left[ \frac{\Gamma (\frac{n+1}4)}{\Gamma (\frac{n+3}4)}
\right]^2\ ,&\text{$n\ge4$ ($k=0$ term).}
\end{cases}
\end{equation}
\end{thm}

\begin{proof}
Using Corollary~\ref{cor3}, determine the maximum of the two terms $k=0,1$ 
in Theorem~\ref{thm:gradient}. 
{From} Corollary~\ref{cor2} for $n=2$ and by explicit calculation for $n=3$, 
one observes that the $k=1$ term is larger. 
For higher dimensions, consider the log of the ratio of
the $k=1$ term to the $k=0$ term: 
set $w= n/4$ and define for $w\ge1$
$$\Lambda (w) = \ln \left\{ \left[\Gamma \Big(w+\frac34\Big)^4
\Big/ \Gamma \Big(w+\frac14\Big)^2 \Gamma \Big(w+\frac54\Big)^2\right] 
\left( 1+\frac4{(4w-1)^2}\right)\right]\ .$$
Observe that by Stirling's formula, $\Lambda (w) \to0$ as $w\to\infty$ and 
$$\Lambda (1) = \ln \left[\frac{117}{25}\ \Gamma\Big(\frac34\Big)^4
\Big/ \Gamma\Big(\frac14\Big)^4\right] 
\simeq - 2.796\ .$$
Since $\Lambda (1)$ is negative, the $k=0$ term is largest for $n=4$.
\begin{equation*}
\begin{split} 
\Lambda'(w) & = 4\left[\psi\Big(w+\frac34\Big) -\psi\Big(w+\frac14\Big)\right]
+\frac{8(4w-1)}{(4w-1)^2+4} - \frac{32w}{16w^2-1}\\
\noalign{\vskip6pt}
& =  2\sum_{k=0}^\infty \left[\Big(k+w+\frac12\Big)^2-\frac1{16}\right]^{-1} 
+ \frac{8(4w-1)}{(4w-1)^2+4} - \frac{32w}{(16w^2-1)}\\
\noalign{\vskip6pt}
& > \frac4{2w+1} + \frac{8(4w-1)}{(4w-1)^2+4} - \frac{32w}{(16w^2-1)} >0
\end{split}
\end{equation*}
for $w\ge1$ since 
$$\frac4{2w+1} + \frac{8(4w-1)}{(4w-1)^2+4} - \frac{32w}{16w^2-1} 
> \frac4{2w+1} + \frac8{4w+1} - \frac8{4w-1} >0$$
for this range of values. 
Hence $\Lambda (w)$ is increasing for $w\ge1$ or $n\ge4$ and since the limit
at infinity is zero, $\Lambda (w)$ must be negative for all $w\ge1$ and 
the $k=0$ term is largest for $n\ge4$. 
This completes the  argument for Theorem~\ref{thm:trace}. 
\renewcommand{\qed}{}
\end{proof}

\begin{thm}\label{thm7}$\quad$

{\rm (A)} For $n-2\le \alpha <n$
$$D_\alpha = \pi^\alpha\left[\Gamma \Big(\frac{n-\alpha}4 +1\Big)\Big/
\Gamma \Big( \frac{n+\alpha}4 +1\Big) \right]^2 
\left( 1+\frac{4\alpha}{(n-\alpha)^2}\right)\ .$$
{\rm (B)} For $n\ge3$ and $\alpha$ sufficiently near $0$ 
$$D_\alpha = \pi^\alpha \left[\Gamma\Big(\frac{n-\alpha}4 +\frac12\Big)
\Big/ \Gamma \Big(\frac{n+\alpha}4 + \frac12\Big)\right]^2\ .$$
Moreover, for $n\ge4$ this value holds when $0<\alpha \le n-3$.

\noindent {\rm (C)}
For large $n$ and fixed $\alpha$, $D_\alpha \simeq (\frac{4\pi}{n})^\alpha$.
\end{thm}

\begin{proof}
The case $n=2$ is contained in Corollary~\ref{cor2}. 
The first step to prove part~(A) will be to set $\alpha = n-2$ and consider 
the log of the ratio of the $k=1$ term to the $k=0$ term. 
Then for 
$$\Lambda = \ln \left\{ \left[\frac{\Gamma (\frac32)\Gamma (\frac{n}2)}
{\Gamma (\frac{n+1}2) \Gamma (1)}\right]^2 (n-1)\right\}\ ,$$
treat $n=w$ as a continuous variable and calculate $\Lambda'(w)$.
\begin{gather*}
\Lambda' (w) = \frac1{w-1} + \psi \Big(\frac{w}2\Big) - \psi\Big(\frac{w+1}2
\Big)\\
\psi\Big(\frac{w+1}2\Big) - \psi\Big(\frac{w}2\Big) = 
2\int_0^\infty \frac{e^{-wt}}{1+e^{-t}}\, dt\ .
\end{gather*}
This formula follows from the Gauss integral representation for $\psi$ (see 
Whittaker and Watson, page~247). 
Set $\delta = w-1$ and write for $\delta >1$
$$\Lambda' = \frac1{\delta} \bigg[ 1-2\delta \int_0^\infty 
\frac{e^{-\delta t}}{1+e^t} \,dt \bigg] 
= \frac2{\delta} \int_0^\infty \frac{e^t e^{-\delta t}}{(1+e^t)^2}\,dt >0\ .$$
Since $\Lambda' >0$, $\Lambda$ as a continuous function of $w$ increases 
from $-\infty$ to $\ln (\pi/2)$. 
$\Lambda (2) =0$ then implies that $\Lambda (n)>0$ for $n\ge3$ and 
verifies the claim in part~(A) for $\alpha =n-2$.
\renewcommand{\qed}{}
\end{proof}

The estimates obtained here for $\psi$ using both Riemann sums and the Gauss
integral are expressed in the following lemma.

\begin{lem}
For $w>1$
\begin{equation}\label{eq17}
\frac2{2w+1} < \psi \Big(\frac{w+1}2\Big) - \psi\Big(\frac{w}2\Big) 
< \begin{cases}
\ds\frac1w + \frac2{w(w+1)}\ ,&\text{$1<w\le 3$}\\
\noalign{\vskip6pt}
\ds\frac1{w-1}\ ,&\text{$3\le w$}\ .
\end{cases}
\end{equation}
\end{lem}

To complete the proof of part (A), set $\alpha = n-2+2\delta$ with 
$0<\delta <1$ and consider 
$$\Lambda (\delta) = \ln \left\{ \left[
\frac{\Gamma [\frac{n+\delta}2] \Gamma [\frac{3-\delta}2]}
{\Gamma [\frac{n+\delta +1}2] \Gamma [1-\frac{\delta}2]} \right]^2
\Big( \frac{n+\delta^2 -1}{(1-\delta)^2}\Big)\right\}\ .$$
Note by using the lemma above 
\begin{gather*}
\Lambda' (\delta) = \psi \Big(\frac{n+\delta}2\Big) 
- \psi\Big(\frac{n+\delta +1}2\Big) +\psi \Big(1-\frac{\delta}2\Big) 
- \psi\Big(\frac32 -\frac{\delta}2\Big) 
+ \frac{2\delta}{n+\delta^2-1} + \frac2{\delta-1}\\
\noalign{\vskip6pt}
>\frac{2\delta}{n+\delta^2-1} + \frac1{1-\delta} - \frac1{n+\delta-1} >0
\end{gather*}
so $\Lambda (\delta)$ is increasing for $0<\delta <1$ and since 
$\Lambda (0)>0$ \ \ $\Lambda(\delta)$ is positive and the result in 
part~(A) is verified.

The purpose of parts (A) and (B) is to demonstrate that for dimension at 
least four or larger there are definite ranges of the parameter $\alpha$ 
where either the $k=0$ or $k=1$ terms give the precise constant for 
Theorem~\ref{thm:gradient}. 
For $n=2$ the $k=1$ term is largest for all $\alpha$. 
For $n=3$, the $k=1$ term is largest except for a small neighborhood of 
$\alpha =0$. 
For $n\ge4$ one can identify a definite interval in the parameter range 
for $\alpha$ where the transition between the two terms for determining 
the maximum value for $D_\alpha$ occurs. 
Consider the log of the ratio of the two terms 
$$\Lambda = \ln \left\{ \left[ 
\frac{\Gamma (\frac{n-\alpha}4 +1)}{\Gamma(\frac{n+\alpha}4 +1)}\ 
\frac{\Gamma (\frac{n+\alpha}4 +\frac12)}{\Gamma(\frac{n-\alpha}4+\frac12)}
\right]^2 
\Big(1+\frac{4\alpha}{(n-\alpha)^2}\Big)\right\}$$
which for $n=3$ becomes 
$$\Lambda = \ln \left\{ \left[
\frac{\Gamma (\frac{7-\alpha}4)\Gamma(\frac{5+\alpha}4)}
{\Gamma (\frac{7+\alpha}4) \Gamma (\frac{5-\alpha}4)} \right]^2 
\Big(\frac{\alpha^2 -2\alpha +9}{(3-\alpha)^2}\Big)\right\}\ .$$
For small values of $\alpha$, $\Lambda$ is still positive which means that 
the $k=1$ term is largest, e.g. for $\alpha = 0.2$, $\Lambda\simeq 0.0021145$,
but for $\alpha = 0,1$, $\Lambda \simeq -0.00103461$. 
Note 
$$\Lambda' (0) = -\psi (7/4) + \psi (5/4) + 4/9 \simeq -0.0304815$$
which requires since $\Lambda (0) =0$ that near zero, $\Lambda (\alpha)<0$ 
and the $k=0$ term is largest. 

To show the case $n\ge4$ for part (B), set $n=\alpha +2\delta$ with 
$\delta \ge 3/2$; then 
$$\Lambda = \ln \left\{ \left[
\frac{\Gamma(\frac{\delta}2+1) \Gamma (\frac{\alpha+\delta}2 +\frac12)}
{\Gamma (\frac{\delta}2 +\frac12) \Gamma (\frac{\alpha+\delta}2 +1)}
\right]^2 \Big( 1+\frac{\alpha}{\delta^2}\Big)\right\}$$
which can now be viewed as a function of $\alpha$ for $0<\alpha \le n-2\delta$
and $\delta \ge 3/2$. 
For this range of values of $\delta$ 
\begin{equation*}
\begin{split}
\Lambda'(\alpha) & = \psi \left(\frac{\alpha+\delta+1}2\right) 
- \psi \left( \frac{\alpha+\delta}2 +1\right) +\frac1{\alpha +\delta^2}\\
\noalign{\vskip6pt}
& < - \frac2{\alpha+\delta +2} + \frac1{\alpha+\delta^2} <0\ .
\end{split}
\end{equation*}
Since $\Lambda (0)=0$, the desired value of $\Lambda$ for $\alpha =n-2\delta$ 
must be negative which will imply that the $k=0$ term is largest for 
$0<\alpha \le n-3$ and $n\ge4$. 
This completes the proof of part~(B).
\medskip

Part (C) follows directly as an application of Stirling's formula 
$$\Gamma (z+a) \simeq \sqrt{2\pi}\ e^{-z} z^{z+a-\frac12} \ \text{ as }\ 
z\to\infty$$
since for fixed $\alpha$ and large  $n$, the value $D_\alpha$ will 
use the $k=0$ term 
$$D_\alpha \simeq \left(\frac{4\pi}{n}\right)^\alpha\ .$$
Note that for $\alpha=2$, this becomes an exact relation for $n\ge5$: 
$D_2 = (4\pi/n)^2$.

\section{Logarithmic uncertainty}

The basic relation for Pitt's inequality with gradient terms 
$$\int_{\real^n} |\nabla f|^2 |x|^{-\alpha} \,dx 
\le 4\pi^2 D_\alpha \int_{\real^n} |\hatf (y)|^2 |y|^{\alpha+2}\,dy$$
becomes an equality at $\alpha=0$ ($D_0 = 1$) so it can be differentiated 
at this value of $\alpha$. 
Corollary~\ref{cor2} and Theorem~\ref{thm7} express the value of $D_\alpha$ 
needed for this calculation.

\begin{thm}\label{thm:log-uncertain}
For $f\in \S(\real^n)$ and $n\ge2$ 
\begin{gather}
\int_{\real^n} \ln |x|\, |\nabla f|^2\,dx 
+ 4\pi^2 \int_{\real^n} \ln |y|\, |y|^2 |\hatf(y)|^2 \,dy 
\ge E\int_{\real^n} |\nabla f|^2\,dx \label{eq18}\\
\noalign{\vskip6pt}
E = \begin{cases} \ds \psi \Big(\frac32\Big) - \ln \pi-1\ ,&\text{$n=2$}\\
\noalign{\vskip6pt}
\ds \psi \Big(\frac{n}4 +\frac12\Big) - \ln\pi\ ,&\text{$n\ge3$\ .}
\end{cases}\notag
\end{gather}
\end{thm}

The increase in the constant here over the corresponding value in 
Theorem~\ref{thm:log} reflects the comparison of the integrals 
$$\int_{\real^n} \ln |x|\, |\nabla f|^2\,dx 
\ge\int_{\real^n} \ln |x|\, \big|\F^{-1} (2\pi |y|\hatf (y))\big|^2\,dx \ .
$$

\section{Iterated Stein-Weiss potentials}

The Stein-Weiss potentials discussed here act at the spectral level, 
that is, in terms of $L^2$ estimates. 
In the context of using these potentials to define a linear operator on 
$L^2(\real^n)$, it is natural to examine iterated applications. 
For example, the Stein-Weiss potential with $0<\alpha <n$
$$\int_{\real^n\times\real^n} \mkern-36mu
f(x) |x|^{-\alpha/2} |x-y|^{-(n-\alpha)} |y|^{-\alpha/2} f(y)\,dx\,dy$$
corresponds to the linear operator 
$$g\too |x|^{-\alpha/2}  (|x|^{-(n-\frac{\alpha}2)} * g)$$
and inequality \eqref{eq:pittsproof} can be rephrased
\begin{equation}\label{eq:pittsproof2} 
\Big\| |x|^{-\alpha/2} (|x|^{-(n-\frac{\alpha}2)} * g)\Big\|_{L^2(\real^n)}
\le \pi^{n/2} 
\left[ \frac{\Gamma (\frac{\alpha}4) \Gamma (\frac{n-\alpha}4)}
{\Gamma (\frac{n}2 - \frac{\alpha}4) \Gamma (\frac{n+\alpha}4)}
\right] \|g\|_{L^2(\real^n)}
\end{equation}
or as a weighted Sobolev inequality
\begin{equation}\label{eq:weightedS}
\|h\|_{L^2(\real^n} \le 
2^{-\alpha/2} \left[ 
\frac{\Gamma (\frac{n-\alpha}4)}{\Gamma(\frac{n+\alpha}4)}\right] 
\|(-\Delta )^{\alpha/4} (|x|^{\alpha/2} h)\|_{L^2(\real^n)}
\end{equation}
which further implies by using the $\|T^*T\| = \|T\|^2$ argument that 
\begin{equation}\label{eq:21}
\|h\|_{L^2(\real^n} \le 2^{-\alpha}
\left[ \frac{\Gamma (\frac{n-\alpha}4)}{\Gamma(\frac{n+\alpha}4)}\right]^2
\big\|\, |x|^{\alpha/2} 
(-\Delta )^{\alpha/2} (|x|^{\alpha/2} h)\big\|_{L^2(\real^n)} 
\quad\footnote{This inequality corresponds to the case $\mu=\lambda$ for 
inequality (7) in \cite{7}.}
\end{equation}

These latter inequalities extend to include successive applications of powers 
of $|x|$ and $(-\Delta)^{1/2}$ and correspond to iterated Stein-Weiss 
potentials subject to growth constraints on the size of the powers and 
the dilation constraint that the sum of the powers of $|x|$ must equal the 
sum of the powers of $(-\Delta)^{1/2}$. 
At the second iteration level, this algorithm leads to a result that 
includes the Maz'ya-Eilertsen inequality:
\begin{equation}\label{eq22}
\| h\|_{L^2(\real^n)} \le 
C\|(-\Delta)^{\rho/4} |x|^{\sigma/2} (-\Delta)^{\beta/4} 
(|x|^{\alpha/2}h) \|_{L^2(\real^n)}
\end{equation}
with $\sigma+\alpha = \rho+\beta$, and 
\addtocounter{footnote}{1}
\footnotetext{The case $\rho=0$, $\sigma=\mu$, $\beta=2\lambda$ and 
$\alpha = 2\lambda-\mu$ corresponds to inequality (7) in \cite{7}.}
\begin{equation}\label{eq23} 
\begin{split}
&\int_{(\real^n)^4}\mkern-18mu 
g(w) |w|^{-\rho/2} |x-w|^{-(n-\sigma/2)} |x|^{-\beta/2} 
|x-y|^{-(n-\alpha)} |y|^{-\beta/2} |y-v|^{-(n-\sigma/2)} |v|^{-\rho/2} \\
&\hskip2truein 
g(v)\, dx\,dy\,dw\,dv \le C_1 \int_{\real^n} |g|^2\,dx\ .\ {\ }^2
\end{split}
\end{equation}
By applying symmetrization to inequality~\eqref{eq23}, one reduces the 
calculation of a sharp constant to considering non-negative radial 
decreasing functions. 
Since $g$ can now be taken to be radial, set $|w|=p$ $|x|=t$, $|y|=s$, 
$|v| =r$ with $u(p) = |w|^{n/2} g(w)$, and inequality~\eqref{eq23} is 
reduced to a convolution inequality on the multiplicative group $\real_+$ 
$$\|\varphi *u\|_{L^2(\real_+)} 
\le \|\varphi\|_{L^1(\real_+)} \|u\|_{L^2(\real_+)}$$
where $\varphi = \kappa *\psi_\alpha *\kappa$ with 
$$\psi_\alpha (t) = \int_{S^{n-1}} \Big[ t+\frac1t -2\xi_1\Big]^{-(n-\alpha)/2}
\,d\xi\ ,\qquad 
\kappa (t) = t^{-\rho/2\ +\ \sigma/4} \psi_{\sigma/2} (t)$$
where $d\xi = $ normalized surface measure on $S^{n-1}$. 
$$\|\varphi\|_{L^1(\real_+)} 
= \|\kappa *\psi_\alpha *\kappa\|_{L^1(\real_+)} 
= (\|\kappa\|_{L^1(\real_+)})^2 \|\psi_\alpha\|_{L^1(\real_+)}$$
$$\|\psi_\alpha\|_{L^1(\real_+)} 
= \frac{\Gamma (\frac{n}2) \Gamma (\frac{\alpha}2)}
{2\Gamma (\frac{n-\alpha}2)} 
\left[ \frac{\Gamma (\frac{n-\alpha}4)}{\Gamma (\frac{n+\alpha}4)}\right]^2\ ,
\qquad 0<\alpha <n$$
\begin{equation*}
\begin{split}
\|\kappa\|_{L^1(\real_+)} 
& = \int_0^\infty t^{-\rho/2\ +\ \sigma/4}
\left[ \int_{S^{n-1}} \Big[ t+\frac1t-2\xi_1\Big]^{-(n-\frac{\sigma}2)/2}
d\xi\right] \frac{dt}t\\
\noalign{\vskip6pt}
&= \left[ \frac{2\pi^{n/2}}{\Gamma (n/2)}\right]^{-1} 
\int_{\real^n} |x-y|^{-(n-\frac{\sigma}2)}  
|y|^{-(\rho+n)/2}dy\ ,\qquad |x|=1\\
\noalign{\vskip6pt}
&= \frac{\Gamma (\frac{n}2) \Gamma (\frac{\sigma}4) \Gamma (\frac{n-\rho}2)
\Gamma [\frac{n+\rho-\sigma}4]}
{2\Gamma (\frac{n}2 -\frac{\sigma}4) \Gamma (\frac{\rho+n}2) 
\Gamma [\frac{n+\sigma-\rho}4]}
\end{split}
\end{equation*}

\begin{equation*}
\begin{split}
C_1 & = \left[ \frac{2\pi^{n/2}}{\Gamma (n/2)}\right]^3 
\|\varphi\|_{L^1(\real_+)} \\
\noalign{\vskip6pt}
& = \pi^{3n/2} \frac{\Gamma(\frac{\alpha}2)}{\Gamma(\frac{n-\alpha}2)}
\left[
\frac{\Gamma (\frac{n-\alpha}4)}{\Gamma (\frac{n+\alpha}4)} 
\frac{\Gamma (\frac{\sigma}4)}{\Gamma (\frac{n}2-\frac{\sigma}4)} 
\frac{\Gamma (\frac{n-\rho}2)}{\Gamma (\frac{n+\rho}2)} 
\frac{\Gamma (\frac{n+\rho-\sigma}4)}{\Gamma (\frac{n+\sigma-\rho}4)} 
\right]^2
\end{split}
\end{equation*}
and 
$$C = 2^{-(\alpha+\sigma)/2}  
\left[
\frac{\Gamma (\frac{n-\alpha}4)}{\Gamma (\frac{n+\alpha}4)} 
\frac{\Gamma (\frac{n-\rho}2)}{\Gamma (\frac{n+\rho}2)} 
\frac{\Gamma (\frac{n+\rho-\sigma}4)}{\Gamma (\frac{n+\sigma-\rho}4)} 
\right]\ .$$
These calculations comprise the proof of the following theorem:

\begin{thm}\label{thm9}
For $g,h\in \S(\real^n)$ and $0<\alpha,\beta,\rho,\sigma <n$, 
$\alpha+\sigma = \beta+\rho$
\begin{gather}
\int_{(\real^n)^4} \mkern-30mu g(w) 
|w|^{-\frac{\rho}2} |x\!-\!w|^{-(n-\frac{\sigma}2)} |x|^{-\frac{\beta}2} 
|x\!-\!y|^{-(n-\alpha)} |y|^{-\frac{\beta}2} |y\!-\!v|^{-(n-\frac{\sigma}2)} 
|v|^{-\frac{\rho}2} g(v)\,dx\,dy\,dw\,dv  \label{eq24}\\
\hskip1.5truein \le B_{\alpha,\rho,\sigma} \int_{\real^n} |g|^2\,dx
\notag\\
\noalign{\vskip6pt}
B_{\alpha,\rho,\sigma} = \pi^{3n/2} 
\frac{\Gamma (\frac{\alpha}2)}{\Gamma (\frac{n-\alpha}2)}
\left[
\frac{\Gamma (\frac{n-\alpha}4)}{\Gamma (\frac{n+\alpha}4)} 
\frac{\Gamma (\frac{\sigma}4)}{\Gamma (\frac{n}2-\frac{\sigma}4)} 
\frac{\Gamma (\frac{n-\rho}2)}{\Gamma (\frac{n+\rho}2)} 
\frac{\Gamma (\frac{n+\rho-\sigma}4)}{\Gamma (\frac{n+\sigma-\rho}4)} 
\right]^2
\notag\\
\noalign{\vskip6pt}
\|h\|_{L^2(\real^n)} \le C_{\alpha,\rho,\sigma} 
\|(-\Delta)^{\rho/4} |x|^{\sigma/2} (-\Delta)^{\beta/4}  (|x|^{\alpha/2}h)
\|_{L^2(\real^n)} \label{eq25}\\
\noalign{\vskip6pt}
C_{\alpha,\rho,\sigma}= 2^{-(\alpha+\sigma)/2} 
\left[
\frac{\Gamma (\frac{n-\alpha}4)}{\Gamma (\frac{n+\alpha}4)} 
\frac{\Gamma (\frac{n-\rho}2)}{\Gamma (\frac{n+\rho}2)} 
\frac{\Gamma (\frac{n+\rho-\sigma}4)}{\Gamma (\frac{n+\sigma-\rho}4)} 
\right]\ .\notag
\end{gather}
\end{thm}

\medskip
\begin{rem}
The papers \cite7 and \cite{15} consider a broader set of problems which 
allow parameters that lie outside the normal range of values for 
fractional integrals.
\end{rem}

\section{Pitt's inequality with iterated gradients}

Pitt's inequality with gradient terms 
$$\int_{\real^n} \Phi (1/|x|) |\nabla f|^2\,dx 
\le 4\pi^2 D_\Phi \int_{\real^n} \Phi (|y|) |y|^2|\hatf (y)|^2\,dy$$
is an intrinsic refinement of the classical inequality and extends 
naturally to iterated gradients 
\begin{equation}\label{eq26} 
\int_{\real^n} \Phi (1/|x|) |\nabla^\ell f|^2\,dx 
\le (4\pi^2)^\ell D_{\Phi,\ell} \int_{\real^n} \Phi (y) |y|^{2\ell} 
|\hatf (y)|^2\,dy
\end{equation}
where 
$$|\nabla^\ell f|^2 = \sum_{p_1=1}^n \cdots \sum_{p_\ell=1}^n 
\left(\frac{\partial}{\partial x_{p_1}} \cdots 
\frac{\partial}{\partial x_{p_\ell}} f\right)^2$$
and $D_{\Phi,\ell} <C_\Phi$ (see equation~\eqref{eq:spectral}).

\begin{thm}\label{thm10}
For $f\in \S(\real^n)$ and $0<\alpha <n$, $n>1$ 
\begin{equation}\label{eq27}
\int_{\real^n} |\nabla^\ell f|^2 |x|^{-\alpha}\,dx 
\le (4\pi^2)^\ell D_{\alpha,\ell} \int_{\real^n} |\hatf (y)|^2 
|y|^{\alpha+2\ell} \,dy\ .
\end{equation}
For $f\in L^2 (\real^n)$
\begin{gather}
\Big| \int_{\real^n\times\real^n}\mkern-18mu 
f(x) \frac1{|x|^{\alpha/2}} \left(\frac{x\cdot y}{|x|\, |y|}\right)^\ell 
\frac1{|x-y|^{n-\alpha}} \frac1{|y|^{\alpha/2}} f(y)\,dx\,dy\Big|
\label{eq28}\\
\noalign{\vskip6pt}
\le \left[ \pi^{\frac{n}2 -\alpha} \Gamma \Big(\frac{\alpha}2\Big)\Big/
\Gamma\Big(\frac{n-\alpha}2\Big)\right] 
D_{\alpha,\ell} \int_{\real^n} |f|^2\,dx\ .\notag
\end{gather}
\end{thm}
\medskip

Observe that when $\ell$ is even, the kernel in \eqref{eq28} is 
positive and the issue of calculating $D_{\alpha,\ell}$ is reduced to 
considering  radial functions.
\medskip

\begin{lem}
For $F,G\in L^2 (S^n)$ and $K (\xi\cdot\eta)\ge 0$ with $\xi,\eta\in S^n$
$$\Big| \int_{S^n\times S^n}F(\xi) K(\xi\cdot\eta) G(\eta)\,d\xi\,d\eta\Big|
\le \bigg(\int_{S^n} K(\xi_1)\,d\xi\bigg) \|F\|_{L^2 (S^n)} \|G\|_{L^2(S^n)}$$
\end{lem}

\begin{proof} 
Split the integrand into the product  of two parts, $F\sqrt{K}$ and 
$G\sqrt{K}$, and apply H\"older's inequality.

For radial functions, inequality~\eqref{eq28} is equivalent to the 
convolution inequality on the group $\real_+$
\begin{gather}
\Big| \int_{\real_+\times\real_+} \mkern-18mu 
h(t) \psi_{\alpha,\ell} (s/t) h(s)\frac{ds}{s}\ \frac{dt}t\Big| 
\le D_{\alpha,\ell} \left[
\frac{\Gamma(\frac{n}2) \Gamma (\frac{\alpha}2)}
{2\pi^\alpha \Gamma (\frac{n-\alpha}2)}\right] 
\int_{\real_+} |h|^2 \frac{dt}t \label{eq29}\\
\noalign{\vskip6pt}
\psi_{\alpha,\ell} (t) = \int_{S^{n-1}} (\xi_1)^\ell 
\Big[ t+\frac15 - 2\xi_1\Big]^{-(n-\alpha)/2}d\xi\notag
\end{gather}
Hence, for $\ell$ even
$$D_{\alpha,\ell} = \frac{2\pi^\alpha \Gamma(\frac{n-\alpha}2)}
{\Gamma (\frac{n}2) \Gamma (\frac{\alpha}2)} 
\|\psi_{\alpha,\ell}\|_{L^1(\real_+)}\ .$$
For $\ell=2$, observe that 
$$\psi_{\alpha,n,2} = \psi_{\alpha,n,0} - \Big(\frac{n-1}n\Big) 
\psi_{\alpha+2,\, n+2,\, 0}$$
and that from equation~\eqref{eq:pitt}
$$D_{\alpha,0} = \pi^\alpha \left[ \frac{\Gamma (\frac{n-\alpha}4)}
{\Gamma (\frac{n+\alpha}4)}\right]^2$$
which then provides: $D_{\alpha,n,2} = D_{\alpha,n,0} - 
\frac{\alpha}4 (n-1) 
D_{\alpha+2,\, n+2,\, 0}$
\renewcommand{\qed}{}
\end{proof}

\begin{Cor}
\begin{equation}\label{eq30} 
D_{\alpha,2} = \pi^\alpha 
\left[ \frac{\Gamma (\frac{n-\alpha}4)}{\Gamma (\frac{n+\alpha}4)}\right]^2
\left[ \frac{(n-\alpha)^2 +4\alpha}{(n+\alpha)^2}\right]\ .
\end{equation}
\end{Cor}
\medskip

More generally, this argument determines a recursion formula for computing 
constants for the case of radial functions which includes the case when 
$\ell$ is an even integer: 
$$D_{\alpha,n,\ell+2} = D_{\alpha,n,\ell} - 
\frac{\alpha}4 (n-1) 
D_{\alpha+2,\, n+2,\, \ell}$$

\section*{Acknowledgements}

The computer 
program {\em mathematica\/} was used to aid some numerical calculations.


\end{document}